# Lyapunov-Based Boundary Feedback Design For Parabolic PDEs

**Iasson Karafyllis**

Dept. of Mathematics, National Technical University of Athens,
Zografou Campus, 15780, Athens, Greece, email: iasonkar@central.ntua.gr

**Abstract**

This paper presents a methodology for the construction of simple Control Lyapunov Functionals (CLFs) for boundary controlled parabolic Partial Differential Equations (PDEs). The proposed methodology provides functionals that contain only simple (and not double or triple) integrals of the state. Moreover, the constructed CLF is "almost diagonal" in the sense that it contains only a finite number of cross-products of the (generalized) Fourier coefficients of the state. The methodology for the construction of a CLF is combined with a novel methodology for boundary feedback design in parabolic PDEs. The proposed feedback design methodology is Lyapunov-based and the feedback controller is an "integral" controller with internal dynamics. It is also shown that the obtained simple CLFs can provide nonlinear boundary feedback laws which achieve global exponential stabilization of semilinear parabolic PDEs with nonlinearities that satisfy a linear growth condition.

## 1. Introduction

One of the most important methodologies for the construction of stabilizing feedback laws in systems described by Ordinary Differential Equations (ODEs) is the Control Lyapunov Function (CLF) methodology (see [23,8] and references therein). The CLF methodology allows the solution of global feedback stabilization problems for highly nonlinear systems and also allows the robustness analysis of the resulting closed-loop system. When studying systems described by Partial Differential Equations (PDEs), the Control Lyapunov Function becomes a Control Lyapunov Functional (CLF). The use of CLFs for the solution of global feedback stabilization problems for systems with PDEs has been presented in detail in [5] and has been used for instance in [1,7,12,14,20,21].

However, it is true that the CLF methodology has not been used so far for systems described by parabolic PDEs with boundary control. Indeed, the design of global boundary feedback stabilizers for parabolic PDEs has focused on linear parabolic PDEs, where several methodologies are available; see [2,6] and the backstepping design methodology described in [13,22]. The backstepping design methodology is essentially a transformation methodology, which can ultimately provide a CLF for the parabolic PDE. However, the CLF depends on the kernel of the transformation (i.e., the computation of the CLF requires the solution of a PDE) and is very complicated even for simple cases: it contains double and triple integrals and it is a very "non-diagonal" functional, in the sense that it contains an infinite number of cross-products of the (generalized) Fourier coefficients of the state. Perhaps, this is the reason that CLFs for boundary controlled parabolic PDEs have not been used so far for feedback design. The complicated nature of the CLF obtained by backstepping has not allowed the use of the CLF methodology for nonlinear parabolic PDEs. Very few feedback design methodologies have been proposed for unstable nonlinear parabolic PDEs: see the extension of the backstepping boundary feedback design in [24,25] as well as feedback designs for distributed inputs in [4,7,20,21]. In many cases the stabilization results are local, guaranteeing exponential stability in specific spatial norms. The papers [10,11] presented methodologies for global feedback stabilization of boundary controlled nonlinear parabolic PDEs: a small-gain methodology is applied in [10], while a CLF methodology is used for PDEs with at most one unstable mode in [11].

In this paper we show how we can obtain simple CLFs for boundary controlled parabolic PDEs. The proposed methodology for the construction of the CLF provides functionals that contain only simple (and not double or triple) integrals of the state. Moreover, the constructed CLF is "almost diagonal" in the sense that it contains only few cross-products of the (generalized) Fourier coefficients of the state. The methodology does not rely on the solution of a PDE but requires knowledge of the eigenvalues and eigenfunctions of a specific Sturm-Liouville (SL) operator.

In order to be able to construct simple CLFs for boundary controlled linear parabolic PDEs, we also present a novel methodology for boundary feedback design. The proposed methodology is Lyapunov-based and is very different from the backstepping methodology. The methodology provides a family of stabilizing boundary feedback laws that have never been used so far: the feedback controller is an "integral" controller with internal dynamics (described by additional ODEs).

As expected the construction of simple CLFs allows the study of global feedback stabilization problems for boundary controlled nonlinear parabolic PDEs. Indeed, we show that simple CLFs can provide nonlinear boundary feedback laws which achieve global exponential stabilization in the $L^2$ norm of semilinear parabolic PDEs with nonlinearities that satisfy a linear growth condition.

The structure of the paper is as follows. Section 2 is devoted to the presentation of the methodology for the construction of the CLF for boundary controlled linear parabolic PDEs. As explained above the methodology is strongly related to a novel methodology for boundary feedback design, which is also presented in detail in Section 2. Section 3 of the present work is devoted to the study of global feedback stabilization problems for boundary controlled semilinear parabolic PDEs with nonlinearities that satisfy a linear growth condition. The proofs of the main results are given in Section 4 of the paper. Finally, the concluding remarks of the paper are provided in Section 5.

**Notation.** Throughout this paper, we adopt the following notation.

* $\Re_+ := [0,+\infty)$. $|y|$ denotes the Euclidean norm of a vector $y \in \Re^n$. $B^T \in \Re^{n \times m}$ denotes the transpose of the matrix $B \in \Re^{m \times n}$. $diag(a_1,...,a_n)$ denotes the diagonal $n \times n$ matrix with the real numbers $a_1,...,a_n$ on its diagonal.

* Let $A \subseteq \Re^n$ be an open set and let $\Omega \subseteq \Re$ and $A \subseteq U \subseteq \bar{A}$ be given sets. By $C^0(U)$ (or $C^0(U;\Omega)$), we denote the class of continuous mappings on $U$ (which take values in $\Omega$). By $C^k(U)$ (or $C^k(U;\Omega)$), where $k \geq 1$, we denote the class of continuous functions on $U$, which have continuous derivatives of order $k$ on $U$ (and also take values in $\Omega$). For a differentiable function $u \in C^0([0,1])$, $u'(x)$ denotes the derivative with respect to $x \in [0,1]$.

* Let $r \in C^0([0,1];(0,+\infty))$ be a positive function. $L^2_r(0,1)$ denotes the Hilbert space of measurable functions $u:[0,1] \to \Re$ with inner product $\langle u,v \rangle = \int_0^1 r(x)u(x)v(x)dx$ for $u,v \in L^2_r(0,1)$. When $r(x) \equiv 1$ then we simply write $L^2(0,1)$. For an interval $I \subseteq \Re_+$, $C^0(I;L^2_r(0,1))$ ($C^1(I;L^2_r(0,1))$) denotes the space of continuous (continuously differentiable) mappings $u: I \to L^2_r(0,1)$.

* Let $u: \Re_+ \times [0,1] \to \Re$ be given. We use the notation $u[t]$ to denote the profile at certain $t \geq 0$, i.e., $(u[t])(x) = u(t,x)$ for all $x \in [0,1]$. When $u(t,x)$ is (twice) differentiable with respect to $x \in [0,1]$, we use the notation $u_x(t,x)$ ($u_{xx}(t,x)$) for the (second) derivative of $u$ with respect to $x \in [0,1]$, i.e., $u_x(t,x) = \frac{\partial u}{\partial x}(t,x)$ ($u_{xx}(t,x) = \frac{\partial^2 u}{\partial x^2}(t,x)$). When $u(t,x)$ is differentiable with respect to $t$, we use the notation $u_t(t,x)$ for the derivative of $u$ with respect to $t$, i.e., $u_t(t,x) = \frac{\partial u}{\partial t}(t,x)$.

* For an integer $k \geq 1$, $H^k(0,1)$ denotes the Sobolev space of functions in $L^2(0,1)$ with all its weak derivatives up to order $k \geq 1$ in $L^2(0,1)$.

## 2. Construction of the Lyapunov Functional

In this section we present a procedure for the construction of a simple Control Lyapunov Functional (CLF) for a 1-D linear parabolic PDE under boundary control. We use the word "simple" in order to indicate that the constructed CLF does not involve any solution of an additional PDE and does not contain double and triple integrals. The only thing that is required for the construction of the CLF is the knowledge of the eigenfunctions of a specific Sturm-Liouville (SL) operator. All steps in the construction are explained in detail.

Consider the SL operator $A: D \to L_r^2(0,1)$ defined by

$$(Af)(x) = -\frac{1}{r(x)} \frac{d}{dx}\left(p(x) \frac{df}{dx}(x)\right) + \frac{q(x)}{r(x)} f(x), \text{ for all } f \in D \text{ and } x \in (0,1) \quad (2.1)$$

where $p \in C^1([0,1]; (0,+\infty))$, $r \in C^0([0,1]; (0,+\infty))$, $q \in C^0([0,1]; \Re)$ and $D \subseteq H^2(0,1)$ is the set of all functions $f: [0,1] \to \Re$ for which

$$b_1 f(0) + b_2 f'(0) = a_1 f(1) + a_2 f'(1) = 0 \quad (2.2)$$

where $a_1, a_2, b_1, b_2$ are real constants with

$$a_1^2 + a_2^2 = 1, \quad b_1^2 + b_2^2 = 1 \quad (2.3)$$

It is well-known (Chapter 11 in [3] and pages 498-505 in [18]) that all eigenvalues of the SL operator $A: D \to L_r^2(0,1)$, defined by (2.1), (2.2) are real. The eigenvalues form an infinite, increasing sequence $\lambda_1 < \lambda_2 < \ldots < \lambda_n < \ldots$ with $\lim_{n \to \infty}(\lambda_n) = +\infty$ and to each eigenvalue $\lambda_n \in \Re$ ($n = 1, 2, \ldots$) corresponds exactly one eigenfunction $\phi_n \in C^2([0,1])$ that satisfies $A\phi_n = \lambda_n \phi_n$, $\|\phi_n\| = 1$ and $b_1 \phi_n(0) + b_2 \phi_n'(0) = a_1 \phi_n(1) + a_2 \phi_n'(1) = 0$. Moreover, the eigenfunctions form an orthonormal basis of $L_r^2(0,1)$.

In the present work, we use the following assumption for the SL operator $A: D \to L_r^2(0,1)$ defined by (2.1), (2.2), (2.3).

**(H):** *The SL operator $A: D \to L_r^2(0,1)$ defined by (2.1), (2.2), (2.3), satisfies*

$$\sum_{n=N+1}^{\infty} \lambda_n^{-1} \max_{0 \leq x \leq 1}(|\varphi_n(x)|) < +\infty, \text{ for certain } N > 0 \text{ with } \lambda_{N+1} > 0 \quad (2.4)$$

It is important to notice that the validity of Assumption (H) can be verified without knowledge of eigenvalues and the eigenfunctions of the SL operator $A$ (see [35]).

In this section, we consider the following control system

$$\frac{\partial u}{\partial t}(t,x) - \frac{1}{r(x)} \frac{\partial}{\partial x}\left(p(x) \frac{\partial u}{\partial x}(t,x)\right) + \frac{q(x)}{r(x)} u(t,x) = 0, \quad x \in (0,1) \quad (2.5)$$

$$b_1 u(t,0) + b_2 \frac{\partial u}{\partial x}(t,0) = a_1 u(t,1) + a_2 \frac{\partial u}{\partial x}(t,1) - \sum_{i=1}^{j} y_i(t) = 0, \tag{2.6}$$

$$\dot{y}_i(t) = \bar{v}_i(t), \quad i = 1,\ldots, j \tag{2.7}$$

where $(u[t], y_1(t),\ldots, y_j(t)) \in L_r^2(0,1) \times \Re^j$ is the state and $(\bar{v}_1(t),\ldots,\bar{v}_j(t)) \in \Re^j$ is the control input. It is clear that system (2.5), (2.6), (2.7) is a 1-D linear parabolic PDE with boundary control given by (2.7) and the relations

$$a_1 u(t,1) + a_2 \frac{\partial u}{\partial x}(t,1) = U(t)$$

$$U(t) = \sum_{i=1}^{j} y_i(t) \tag{2.8}$$

Let $\mu_i > 0$, $i = 1,\ldots, j$ with $\mu_i \neq \lambda_n$ for $n = 1, 2,\ldots$, $i = 1,\ldots, j$ be given constants and consider functions $\varphi_i \in C^2([0,1])$, $i = 1,\ldots, j$ that satisfy

$$\left(p(x)\varphi_i'(x)\right)' - q(x)\varphi_i(x) = -\mu_i r(x)\varphi_i(x), \quad x \in [0,1], \quad i = 1,\ldots, j \tag{2.9}$$

$$b_1 \varphi_i(0) + b_2 \varphi_i'(0) = a_1 \varphi_i(1) + a_2 \varphi_i'(1) - 1 = 0, \tag{2.10}$$

Since $\mu_i > 0$, $i = 1,\ldots, j$ with $\mu_i \neq \lambda_n$ for $n = 1, 2,\ldots$, $i = 1,\ldots, j$, it follows that each of the boundary-value problems (2.9), (2.10) has a unique solution.

We perform the state transformation

$$w = u - \sum_{i=1}^{j} \varphi_i y_i \tag{2.11}$$

and the input transformation

$$\bar{v}_i = -\mu_i y_i + v_i, \quad i = 1,\ldots, j \tag{2.12}$$

System (2.5), (2.6), (2.7) in the new coordinates is described by the equations

$$\frac{\partial w}{\partial t}(t,x) = -(Aw[t])(x) - \sum_{i=1}^{j} \varphi_i(x) v_i(t), \quad x \in (0,1) \tag{2.13}$$

$$b_1 w(t,0) + b_2 \frac{\partial w}{\partial x}(t,0) = a_1 w(t,1) + a_2 \frac{\partial w}{\partial x}(t,1) = 0, \tag{2.14}$$

$$\dot{y}_i(t) = -\mu_i y_i(t) + v_i(t), \quad i = 1,\ldots, j \tag{2.15}$$

The object of the study of the present section is the PDE-ODE system (2.13), (2.14), (2.15). However, it is clear that the transformation (2.11), (2.12) can always allow the interpretation of the results in the original coordinates (i.e., for the original 1-D linear parabolic system (2.5), (2.6), (2.7)).

Using the fact that $A\phi_n = \lambda_n \phi_n$, we guarantee that the following equations hold for all $n = 1, 2, \ldots$ and $(w, v) \in D \times \Re$:

$$\left\langle \phi_n, -Aw - \sum_{i=1}^{j} v_i \varphi_i \right\rangle = -\lambda_n \langle \phi_n, w \rangle - \sum_{i=1}^{j} v_i \langle \varphi_i, \phi_n \rangle \tag{2.16}$$

Therefore we are led to the study of the following linear, time-invariant, finite-dimensional system

$$\dot{z} = Cz + \sum_{i=1}^{j} B_i v_i \qquad (2.17)$$

$$z \in \mathfrak{R}^N, v = (v_1, ..., v_j) \in \mathfrak{R}^j$$

where $N > 0$ is an integer with $\lambda_{N+1} > 0$ and

$$C = -diag(\lambda_1, ..., \lambda_N), \quad B_i = -\begin{bmatrix} \langle \varphi_i, \phi_1 \rangle \\ \vdots \\ \langle \varphi_i, \phi_N \rangle \end{bmatrix}, \quad i = 1, ..., j \qquad (2.18)$$

For the linear time-invariant system (2.17), we are in a position to establish the following result.

**Lemma 2.1:** *System (2.17) is controllable.*

Since system (2.17) is controllable, there exist vectors $K_i = (K_{i,1}, ..., K_{i,N})^T \in \mathfrak{R}^N$, a constant $\sigma > 0$ and a symmetric, positive definite matrix $R = \{R_{n,m} : n, m = 1, ..., N\} \in \mathfrak{R}^{N \times N}$ such that the following matrix inequality holds:

$$R\left(C + \sum_{i=1}^{j} B_i K_i^T\right) + \left(C + \sum_{i=1}^{j} B_i K_i^T\right)^T R \leq -2\sigma I \qquad (2.19)$$

Therefore, there exist constants $c_2 \geq c_1 > 0$ such that the following inequality holds for all $w \in L_r^2(0,1)$:

$$c_1 \sum_{n=1}^{N} \langle \phi_n, w \rangle^2 \leq \sum_{n=1}^{N} \sum_{m=1}^{N} R_{n,m} \langle \phi_n, w \rangle \langle \phi_m, w \rangle \leq c_2 \sum_{n=1}^{N} \langle \phi_n, w \rangle^2 \qquad (2.20)$$

Next define the linear, continuous operator $G: L^2(0,1) \to L^2(0,1)$ by means of the formula:

$$Gw = \sum_{n=1}^{N} \sum_{m=1}^{N} R_{n,m} \phi_m \langle \phi_n, w \rangle, \text{ for all } w \in L_r^2(0,1) \qquad (2.21)$$

The fact that $R = \{R_{n,m} : n, m = 1, ..., N\} \in \mathfrak{R}^{N \times N}$ is symmetric implies that

$$\langle Gw, u \rangle = \langle Gu, w \rangle, \text{ for all } w, u \in L_r^2(0,1) \qquad (2.22)$$

Moreover, inequalities (2.19), (2.20) in conjunction with definitions (2.18), (2.21) imply that the following inequalities hold for all $(w, v) \in D \times \mathfrak{R}$:

$$c_1 \sum_{n=1}^{N} \langle \phi_n, w \rangle^2 \leq \langle Gw, w \rangle \leq c_2 \sum_{n=1}^{N} \langle \phi_n, w \rangle^2 \qquad (2.23)$$

$$\left\langle Gw, -Aw - \sum_{i=1}^{j} \varphi_i v_i \right\rangle \leq -\sigma \sum_{n=1}^{N} \langle \phi_n, w \rangle^2 - \sum_{i=1}^{j} \left( v_i - \sum_{n=1}^{N} K_{i,n} \langle \phi_n, w \rangle \right) \langle Gw, \varphi_i \rangle \qquad (2.24)$$

Define the projection operator:

$$Pw := \sum_{n=1}^{N} \langle \phi_n, w \rangle \phi_n, \text{ for all } w \in L_r^2(0,1) \qquad (2.25)$$

Let $\omega_1, ..., \omega_j, \gamma > 0$ be constants (to be selected). We define the functional

$$V(w, y) = \frac{1}{2} \langle Gw, w \rangle + \frac{\gamma}{2} \|w - Pw\|^2 + \frac{1}{2} \sum_{i=1}^{j} \omega_i y_i^2, \text{ for all } (w, y) \in L_r^2(0,1) \times \mathfrak{R}^j \qquad (2.26)$$

where $y = (y_1,...,y_j)^T \in \Re^j$. It is clear that inequalities (2.23) in conjunction with equations (2.25), (2.26) imply the inequalities:

$$\frac{\min(c_2, \gamma, \omega_1, ..., \omega_j)}{2}\left(\|w\|^2 + |y|^2\right) \leq V(w, y) \leq \frac{\max(c_2, \gamma, \omega_1, ..., \omega_j)}{2}\left(\|w\|^2 + |y|^2\right), \text{ for all } (w, y) \in L_r^2(0,1) \times \Re^j \quad (2.27)$$

Moreover, it should be noted that the functional $V(w, y)$ defined by (2.26) involves only single integrals (and not double or triple integrals). This becomes apparent by equations (2.21), (2.26) and the fact that $\|w - Pw\|^2 = \|w\|^2 - \sum_{n=1}^{N}\langle \phi_n, w\rangle^2$ (a consequence of Parseval's identity and definition (2.25)), which imply that the following equation holds for all $(w, y) \in L_r^2(0,1) \times \Re^j$:

$$V(w, y) = \frac{1}{2}\sum_{n=1}^{N}\sum_{m=1}^{N} R_{n,m}\langle \phi_n, w\rangle\langle \phi_m, w\rangle + \frac{\gamma}{2}\left(\langle w, w\rangle^2 - \sum_{n=1}^{N}\langle \phi_n, w\rangle^2\right) + \frac{1}{2}\sum_{i=1}^{j}\omega_i y_i^2 \quad (2.28)$$

Using (2.22), (2.26), (2.13), (2.14), (2.15) the derivative $\dot{V}$ of $V(w, y)$ along the solutions of (2.13), (2.14), (2.15) is given by the formula

$$\dot{V} = \left\langle Gw, -Aw - \sum_{i=1}^{j} v_i \varphi_i \right\rangle + \gamma\left\langle w - Pw, -Aw + PAw - \sum_{i=1}^{j} v_i \varphi_i + \sum_{i=1}^{j} v_i P\varphi_i \right\rangle + \sum_{i=1}^{j}\omega_i y_i v_i - \sum_{i=1}^{j}\mu_i \omega_i y_i^2 \quad (2.29)$$

The following theorem shows that if the constants $\omega_1, ..., \omega_j, \gamma > 0$ are selected in an appropriate way then the functional $V(w, y)$ defined by (2.26) is a CLF for system (2.13), (2.14), (2.15).

**Theorem 2.2:** *Assume that* $\omega_1, ..., \omega_j, \gamma > 0$ *satisfy*

$$\sigma\mu_i \geq 2j\omega_i |K_i|^2, \ i = 1, ..., j \quad (2.30)$$

$$\sigma\lambda_{N+1} \geq 2j\gamma\sum_{i=1}^{j}\|\varphi_i\|^2 |K_i|^2 \quad (2.31)$$

*Then* $V(w, y)$ *as defined by (2.26) is a CLF for system (2.13), (2.14), (2.15), in the sense that for every* $L_i \geq 0$, $i = 1, ..., j$ *and for every integer* $M \geq N+1$ *with* $4(\lambda_{M+1} - \lambda_{N+1}) \geq \gamma\sum_{i=1}^{j}L_i \sum_{n=M+1}^{\infty}\langle \phi_n, \varphi_i\rangle^2$ *the following inequality holds for all* $(w, y) \in D \times \Re^j$:

$$\dot{V} \leq -\frac{1}{2}\sum_{i=1}^{j}\omega_i \mu_i y_i^2 - \frac{\min(\gamma\lambda_{N+1}, \sigma)}{2}\|w\|^2 \quad (2.32)$$

*where* $\dot{V}$ *is given by (2.29) and*

$$v_i = \langle k_i, w\rangle - \omega_i L_i y_i, \ i = 1, ..., j \quad (2.33)$$

*with*

$$k_i := \sum_{n=1}^{N}\left(K_{i,n} + L_i\left\langle \sum_{m=1}^{N} R_{n,m}\phi_m, \varphi_i\right\rangle\right)\phi_n + \gamma L_i \sum_{n=N+1}^{M}\langle \phi_n, \varphi_i\rangle\phi_n, \ i = 1, ..., j \quad (2.34)$$

*Moreover, the feedback law (2.33) is a globally exponentially stabilizing feedback law for (2.13), (2.14), (2.15), in the sense that there exist constants* $\bar{K}, \bar{\sigma} > 0$ *such that for every* $w_0 \in D$, $y_0 \in \Re^j$, *the initial-boundary value problem (2.13), (2.14), (2.15), (2.33) with initial condition* $w[0] = w_0$,

$y(0) = y_0$ *admits a unique solution* $w \in C^0(\mathfrak{R}_+ \times [0,1]) \cap C^1((0,+\infty) \times [0,1])$, $y \in C^1(\mathfrak{R}_+; \mathfrak{R}^j)$ *with* $w[t] \in D \cap C^2([0,1])$ *for all* $t > 0$, *which satisfies (2.13), (2.14), (2.15), (2.33) for all* $t > 0$ *and the following estimate*

$$\|w[t]\| + |y(t)| \leq \bar{K} \exp(-\bar{\sigma} t)(\|w_0\| + |y_0|), \text{ for all } t \geq 0 \qquad (2.35)$$

**Remark 2.3:** If $L_i = 0$, $i = 1,...,j$ then formulas (2.33), (2.34) give us the reduced-model boundary feedback design $v_i = \left\langle \sum_{n=1}^{N} K_{i,n} \phi_n, w \right\rangle$, $i = 1,...,j$ (see [6]). Therefore, formulas (2.33), (2.34) contain the reduced-model boundary feedback design as a subcase. However, formulas (2.33), (2.34) give a parameterized family of stabilizing boundary feedback laws which cannot be obtained by backstepping.

The following example illustrates how easily we can use Theorem 2.2 for the construction of a CLF for a linear parabolic PDE.

**Example 2.4:** Consider the reaction-diffusion PDE with Dirichlet actuation:

$$\begin{aligned} u_t &= p u_{xx} - q u \\ u(t,0) &= 0 \\ u(t,1) &= U(t) \end{aligned} \qquad (2.36)$$

for $(t,x) \in (0,+\infty) \times (0,1)$, where $p > 0$, $q \in \mathfrak{R}$ are constants with $-4p\pi^2 < q < -p\pi^2$. The corresponding SL operator is given by (2.1), (2.2), (2.3) with $p(x) \equiv p$, $q(x) \equiv q$, $r(x) \equiv 1$, $b_1 = a_1 = 1$ and $b_2 = a_2 = 0$. The eigenfunctions of the corresponding SL operator are $\phi_n(x) = \sqrt{2}\sin(n\pi x)$ for $n = 1,2,...$ and the eigenvalues are $\lambda_n = pn^2\pi^2 + q$ for $n = 1,2,...$. Therefore Assumption (H) holds. The assumption $-4p\pi^2 < q < -p\pi^2$ allows us to select $N = j = 1$. By selecting $\varphi(x) = \sin\left(\frac{5\pi}{2}x\right)$ and $\mu = p\frac{25\pi^2}{4} + q$, we guarantee that (2.9), (2.10) holds. Therefore, by means of transformation (2.11), (2.12), we are led to the study of the following systems:

$$\begin{aligned} w_t &= p w_{xx} - q w - v \\ w(t,0) &= w(t,1) = 0 \\ \dot{y} &= -\left(p\frac{25\pi^2}{4} + q\right)y + v \end{aligned} \qquad (2.37)$$

and

$$\begin{aligned} \dot{z} &= -(p\pi^2 + q)z + \frac{4\sqrt{2}}{21\pi} v \\ z &\in \mathfrak{R}, v \in \mathfrak{R} \end{aligned} \qquad (2.38)$$

Clearly, (2.19) holds with arbitrary $\sigma > 0$, $R = [1]$ and $K = -\frac{21\pi}{4\sqrt{2}}(\sigma - p\pi^2 - q)$. It follows from Theorem 2.2 that a family of CLFs for system (2.37) is given by

$$V(w,y) = \frac{1-\gamma}{2}\left(\sqrt{2}\int_0^1 \sin(\pi x) w(x) dx\right)^2 + \frac{\gamma}{2}\int_0^1 w^2(x) dx + \frac{\omega}{2} y^2 \qquad (2.39)$$

provided that the constants $\omega, \gamma > 0$ satisfy the inequalities $4\sigma(25 p\pi^2 + 4q) \geq 441\pi^2 (\sigma - p\pi^2 - q)^2 \omega$ and $128 p\pi^2 \sigma \geq 441\pi^2 \gamma (\sigma - p\pi^2 - q)^2$; these are inequalities (2.30), (2.31). Using (2.11) and transforming back to the original coordinates, we obtain the CLF

$$V(u, y) = \frac{1-\gamma}{2}\left(\sqrt{2}\int_0^1 \sin(\pi x) u(x) dx + \frac{4\sqrt{2}}{21\pi} y\right)^2 + \frac{\gamma}{2}\int_0^1 u^2(x) dx + \frac{2\omega y^2 + \gamma}{4} y^2 - \gamma y \int_0^1 u(x) \sin\left(\frac{5\pi}{2} x\right) dx \quad (2.40)$$

for the following PDE system with boundary control

$$\begin{aligned} u_t &= p u_{xx} - q u \\ u(t,0) &= 0 \\ u(t,1) &= y(t) \\ \dot{y} &= -\left(p \frac{25\pi^2}{4} + q\right) y + v \end{aligned} \quad (2.41)$$

As noticed previously, the CLF given by (2.40) involves no double or triple integrals. Moreover, using formulas (2.11), (2.33), (2.34) we obtain a family of feedback stabilizers for system (2.41); namely the parameterized family of feedback laws

$$v(t) = -\left(\omega L + \int_0^1 k(x) \sin\left(\frac{5\pi x}{2}\right) dx\right) y(t) + \int_0^1 k(x) u(t, x) dx \quad (2.42)$$

with

$$k(x) := -\left(\frac{21\pi}{4}(\sigma - p\pi^2 - q) + \frac{8L}{21\pi}\right) \sin(\pi x) + \frac{8\gamma L}{\pi} \sum_{n=2}^{M} \frac{(-1)^n n}{25 - 4n^2} \sin(n\pi x) \quad (2.43)$$

where $L \geq 0$ and $M \geq 2$ is an integer for which the inequality $p\pi^4 ((M+1)^2 - (N+1)^2) \geq 8\gamma L \sum_{n=M+1}^{\infty} \frac{n^2}{(25 - 4n^2)^2}$ holds. ◁

## 3. Nonlinear Feedback is Better than Linear Feedback

When $j = N$, then the matrix $B = [B_1 \ B_2 \ \ldots \ B_N]$, with $B_i \in \mathfrak{R}^N$ ($i = 1, \ldots, N$) being defined by (2.18), is a square matrix. Moreover, if $\det(B) \neq 0$, then (2.19) holds with any diagonal positive definite matrix $R$ and $K = -B^{-1} \text{diag}(\sigma_1 - \lambda_1, \ldots, \sigma_N - \lambda_N)$, where $K = \begin{bmatrix} K_1^T \\ \vdots \\ K_N^T \end{bmatrix}$ and $\sigma_i > 0$, $i = 1, \ldots, N$. Therefore, in this case we are in a position to assign at will the dynamics of the first $N$ modes of the PDE (2.13) with boundary conditions given by (2.14). Moreover, for $R = I$, it follows from (2.21) and (2.25) that $G = P$ and consequently the CLF defined by (2.26) has a particularly simple form. The simplicity of the CLF may be exploited for the study of the semilinear PDE problem

$$\frac{\partial u}{\partial t}(t, x) + (A u[t])(x) = F(u[t])(x), \quad x \in (0,1) \quad (3.1)$$

$$b_1 u(t, 0) + b_2 \frac{\partial u}{\partial x}(t, 0) = a_1 u(t, 1) + a_2 \frac{\partial u}{\partial x}(t, 1) - \sum_{i=1}^{N} y_i(t) = 0, \quad (3.2)$$

$$\dot{y}_i(t) = \bar{v}_i(t), \quad i = 1,\ldots,N \tag{3.3}$$

where $(u[t], y_1(t),\ldots, y_N(t)) \in L_r^2(0,1) \times \Re^N$ is the state, $(\bar{v}_1(t),\ldots,\bar{v}_N(t)) \in \Re^N$ is the control input and $F: L_r^2(0,1) \to L_r^2(0,1)$ is a continuous mapping with $F(0) = 0$ that satisfies the linear growth condition

$$\|F(u)\| \leq \bar{L}\|u\|, \text{ for all } u \in L_r^2(0,1) \tag{3.4}$$

where $\bar{L} \geq 0$ is a constant.

In this section we assume that we have $N$ functions $\varphi_i \in C^2([0,1])$, $i = 1,\ldots,N$ and constants $\mu_i > 0$, $i = 1,\ldots,N$ so that (2.9), (2.10) hold. Moreover, the matrix $B = [B_1 \; B_2 \; \ldots \; B_N]$ with $B_i \in \Re^N$, $i = 1,\ldots,N$ defined by (2.18) satisfies $\det(B) \neq 0$. The state transformation (2.11) with $j = N$ and the input transformation (2.12) with $j = N$ gives us the equivalent system

$$\frac{\partial w}{\partial t}(t,x) + (Aw[t])(x) = \left(F\left(w[t] + \sum_{i=1}^{N} \varphi_i y_i(t)\right)\right)(x) - \sum_{i=1}^{N} \varphi_i(x) v_i(t), \quad x \in (0,1) \tag{3.5}$$

$$b_1 w(t,0) + b_2 \frac{\partial w}{\partial x}(t,0) = a_1 w(t,1) + a_2 \frac{\partial w}{\partial x}(t,1) = 0, \tag{3.6}$$

$$\dot{y}_i(t) = -\mu_i y_i(t) + v_i(t), \quad i = 1,\ldots,N \tag{3.7}$$

Moreover, we will assume that the $N$ functions $\varphi_i \in C^2([0,1])$, $i = 1,\ldots,N$ are orthogonal to each other, i.e., the following implication holds:

$$i \neq m \Rightarrow \langle \varphi_i, \varphi_m \rangle = 0 \tag{3.8}$$

Assumption (3.9) always holds if the $N$ functions $\varphi_i \in C^2([0,1])$, $i = 1,\ldots,N$ are selected to be eigenfunctions of the SL operator $\bar{A}: \bar{D} \to L_r^2(0,1)$ defined by (2.1) with $\bar{D} \subseteq H^2(0,1)$ being the set of all functions $f: [0,1] \to \Re$ for which $b_1 f(0) + b_2 f'(0) = -a_2 f(1) + a_1 f'(1) = 0$ (notice the difference with (2.2)).

We consider two different controllers for system (3.5), (3.6), (3.7).

1st Controller: A Nonlinear Controller

We consider the nonlinear controller (based on cancellation of the nonlinearities for the first $N$ modes)

$$v_i = \sum_{m=1}^{N} g_{i,m}(\sigma - \lambda_m)\langle \phi_m, w \rangle + \sum_{m=1}^{N} g_{i,m} \langle \phi_m, F(u) \rangle, \quad i = 1,\ldots,N \tag{3.9}$$

where $u = w + \sum_{i=1}^{N} \varphi_i y_i$, $g = \{g_{i,m}: i,m = 1,\ldots,N\} = -B^{-1}$ and $\sigma > 0$. The nonlinear controller (3.9) guarantees that the following equations hold for all $n = 1,\ldots,N$ and $w \in D$:

$$\left\langle \phi_n, -Aw - \sum_{i=1}^{N} v_i \varphi_i + F(u) \right\rangle = -\sigma \langle \phi_n, w \rangle \tag{3.10}$$

where $u = w + \sum_{i=1}^{N} \varphi_i y_i$. Therefore, we conclude from (3.10) that the nonlinear controller achieves feedback linearization of the system that describes the evolution of the first $N$ modes of the solution of (3.5), (3.6), (3.7).

## 2nd Controller: A Linear Controller

We consider the linear controller (based on domination for the first $N$ modes)

$$v_i = \sum_{m=1}^{N} g_{i,m} (\sigma - \lambda_m)\langle \phi_m, w \rangle, \quad i = 1,...,N \tag{3.11}$$

where $g = \{g_{i,m} : i,m = 1,...,N\} = -B^{-1}$ and $\sigma > 0$. The linear controller (3.11) guarantees that the following equations hold for all $n = 1,...,N$ and $w \in D$:

$$\left\langle \phi_n, -Aw - \sum_{i=1}^{N} v_i \phi_i + F(u) \right\rangle = -\sigma \langle \phi_n, w \rangle + \langle \phi_n, F(u) \rangle \tag{3.12}$$

where $u = w + \sum_{i=1}^{N} \phi_i y_i$. Notice that for this controller, the system that describes the evolution of the first $N$ modes of the solution of (3.5), (3.6), (3.7) is nonlinear.

For each one of the two controllers, we have the following results.

**Theorem 3.1:** *Suppose that there exists a constant $\kappa > 0$ such that the following inequalities hold:*

$$\mu_i^2 > 2N\bar{L}^2 \left(1 + \frac{1}{\kappa}\right) \|\phi_i\|^2 \sum_{m=1}^{N} g_{i,m}^2, \quad i = 1,...,N \tag{3.13}$$

$$\lambda_{N+1}^2 > \bar{L}^2 (1 + \kappa N)\left(1 + 2N \sum_{i=1}^{N}\sum_{m=1}^{N} \|\phi_i\|^2 g_{i,m}^2\right) \tag{3.14}$$

*Then there exist constants $\omega_1,...,\omega_N, \gamma, R > 0$ such that the following functional*

$$V(w, y) := \frac{R}{2}\sum_{m=1}^{N}\langle \phi_m, w \rangle^2 + \frac{\gamma}{2}\|w - Pw\|^2 + \frac{1}{2}\sum_{i=1}^{N}\omega_i y_i^2 \tag{3.15}$$

*is a CLF for system (3.5), (3.6), (3.7). More specifically, there exists a constant $\theta > 0$ such that the following inequality holds for all $(w, y) \in D \times \Re^N$:*

$$\dot{V} \leq -\theta \left(\|w\|^2 + \sum_{i=1}^{N} y_i^2\right) \tag{3.16}$$

*where*

$$\dot{V} = R\sum_{n=1}^{N}\langle \phi_n, w\rangle\langle \phi_n, \dot{w}\rangle + \sum_{i=1}^{N}\omega_i y_i v_i - \sum_{i=1}^{N}\mu_i \omega_i y_i^2 + \gamma\langle w - Pw, \dot{w} - P\dot{w}\rangle \tag{3.17}$$

$\dot{w} = -Aw - \sum_{i=1}^{N} v_i \phi_i + F(u)$, $u = w + \sum_{i=1}^{N}\phi_i y_i$ *and the control inputs $v_i$, $i = 1,...,N$ are given by (3.9). Furthermore, there exist constants $\bar{K}, \bar{\sigma} > 0$ such that every solution $w \in C^0(\Re_+; L_r^2(0,1)) \cap C^1((0,+\infty); L_r^2(0,1))$, $y \in C^1(\Re_+; \Re^N)$ with $w[t] \in D$ for all $t > 0$, which satisfies (3.5), (3.6), (3.7), (3.9) for all $t > 0$ also satisfies estimate (2.35).*

**Theorem 3.2:** *Suppose that there exists a constant $\kappa > 0$ such that the following inequalities hold:*

$$\sigma^2 > \bar{L}^2(1+\kappa N) \tag{3.18}$$

$$\lambda_{N+1}^2 > \bar{L}^2(1+\kappa N)\left(1+2N\sum_{i=1}^{N}\sum_{m=1}^{N}\frac{\|\varphi_i\|^2 g_{i,m}^2 (\sigma-\lambda_m)^2}{\sigma^2 - \bar{L}^2(1+\kappa N)}\right) \tag{3.19}$$

$$\mu_i^2 > 2N\bar{L}^2\left(1+\frac{1}{\kappa}\right)\|\varphi_i\|^2 \sum_{m=1}^{N}\frac{g_{i,m}^2(\sigma-\lambda_m)^2}{\sigma^2-\bar{L}^2(1+\kappa N)}, \quad i=1,...,N \tag{3.20}$$

*Then there exist constants $\omega_1,...,\omega_N, \gamma, R > 0$ such that the functional $V(w,y)$ defined by (3.15) is a CLF for system (3.5), (3.6), (3.7). More specifically, there exists a constant $\theta > 0$ such that the inequality (3.16) holds for all $(w,y) \in D \times \Re^N$, where $\dot{V}$ is given by (3.17), $\dot{w} = -Aw - \sum_{i=1}^{N} v_i \varphi_i + F(u)$, $u = w + \sum_{i=1}^{N} \varphi_i y_i$ and the control inputs $v_i$, $i=1,...,N$ are given by (3.11). Furthermore, there exist constants $\bar{K}, \bar{\sigma} > 0$ such that every solution $w \in C^0(\Re_+; L_r^2(0,1)) \cap C^1((0,+\infty); L_r^2(0,1))$, $y \in C^1(\Re_+; \Re^N)$ with $w[t] \in D$ for all $t > 0$, which satisfies (3.5), (3.6), (3.7), (3.11) for all $t > 0$ also satisfies estimate (2.35).*

It is clear that when $\lambda_m < 0$ for $m=1,...,N$ then $\dfrac{(\sigma-\lambda_m)^2}{\sigma^2-\bar{L}^2(1+\kappa N)} > 1$ for all $\sigma > 0$ with $\sigma^2 > \bar{L}^2(1+\kappa N)$ and $m=1,...,N$. Therefore, *no matter what $\sigma > 0$ is,* conditions (3.19), (3.20) are more demanding than conditions (3.13), (3.14). This means that the linear controller (3.11) achieves global stabilization of system (3.5), (3.6), (3.7) for a strictly smaller set of nonlinearities than the set of nonlinearities allowed for global stabilization of system (3.5), (3.6), (3.7) under the nonlinear controller (3.9).

Conditions (3.13), (3.14) are equivalent to the single condition

$$\bar{L}^2 < \frac{2\bar{a}\bar{b}}{\bar{a}+\bar{b}+\sqrt{(\bar{a}-\bar{b})^2+4N\bar{a}\bar{b}}} \tag{3.21}$$

where $\bar{a} := \min_{i=1,...,N}\left(\dfrac{\mu_i^2}{2N\|\varphi_i\|^2 \sum_{m=1}^{N} g_{i,m}^2}\right)$ and $\bar{b} := \dfrac{\lambda_{N+1}^2}{1+2N\sum_{i=1}^{N}\sum_{m=1}^{N}\|\varphi_i\|^2 g_{i,m}^2}$. Formula (3.21) may be used in a straightforward way for the computation of an upper bound on the growth coefficient $\bar{L}$ that is allowed for global stabilization of system (3.5), (3.6), (3.7) under the nonlinear controller (3.9). The following example illustrates how easily condition (3.21) may be used for the computation of an upper bound on the growth coefficient $\bar{L}$.

**Example 3.3:** Consider the reaction-diffusion PDE with Dirichlet actuation:

$$\begin{aligned} u_t &= u_{xx} + 5\pi^2 u + F(u) \\ u(t,0) &= 0 \\ u(t,1) &= U(t) \end{aligned} \tag{3.22}$$

where $F \in C^0(\Re)$ is a function for which the linear growth condition

$$|F(u)| \leq \bar{L}|u|, \text{ for all } u \in \Re \tag{3.23}$$

holds for some constant $\bar{L} \geq 0$. The corresponding SL operator is given by (2.1), (2.2), (2.3) with $p(x) \equiv 1$, $q(x) \equiv -5\pi^2$, $r(x) \equiv 1$, $b_1 = a_1 = 1$ and $b_2 = a_2 = 0$. The eigenfunctions of the corresponding SL operator are $\phi_n(x) = \sqrt{2}\sin(n\pi x)$ for $n = 1, 2, ...$ and the eigenvalues are $\lambda_n = (n^2 - 5)\pi^2$ for $n = 1, 2, ...$. Therefore Assumption (H) holds. It is clear that the linear part of the open-loop system (3.22) has two unstable modes. Since $\lambda_3 = 4\pi^2 > 0$, we select $N = 2$. By selecting $\varphi_1(x) = \sin\left(\frac{5\pi}{2}x\right)$, $\varphi_2(x) = -\sin\left(\frac{7\pi}{2}x\right)$ and $\mu_1 = \frac{5\pi^2}{4}$, $\mu_2 = \frac{29\pi^2}{4}$, we guarantee that (2.9), (2.10) holds.

The state transformation (2.11) with $j = 2$ and the input transformation (2.12) with $j = 2$ gives us the equivalent system

$$\frac{\partial w}{\partial t}(t,x) = \frac{\partial^2 w}{\partial x^2}(t,x) + 5\pi^2 w(t,x) - \varphi_1(x)v_1(t) - \varphi_2(x)v_2(t) + f\left(w(t,x) + \varphi_1(x)y_1(t) + \varphi_2(x)y_2(t)\right), \quad x \in (0,1) \quad (3.24)$$

$$w(t,0) = w(t,1) = 0, \quad (3.25)$$

$$\dot{y}_1(t) = -\frac{5\pi^2}{4}y_1(t) + v_1(t), \quad \dot{y}_2(t) = -\frac{29\pi^2}{4}y_2(t) + v_2(t) \quad (3.26)$$

with

$$U(t) = y_1(t) + y_2(t) \quad (3.27)$$

The matrix $B = [B_1 \ B_2]$, with $B_i \in \mathfrak{R}^2$ ($i = 1, 2$) being defined by (2.18), is

$$B = \frac{4\sqrt{2}}{3\pi}\begin{bmatrix} \frac{1}{7} & \frac{1}{15} \\ -\frac{2}{3} & -\frac{2}{11} \end{bmatrix} \quad (3.28)$$

and $\det(B) \neq 0$ holds. The coefficients $\{g_{i,m} : i, m = 1, 2\}$ involved in (3.9) are

$$g_{1,1} = \frac{1890\pi}{256\sqrt{2}}, \quad g_{1,2} = \frac{693\pi}{256\sqrt{2}}, \quad g_{2,1} = -\frac{6930\pi}{256\sqrt{2}}, \quad g_{2,2} = -\frac{1485\pi}{256\sqrt{2}} \quad (3.29)$$

Using (3.21) we are in a position to guarantee that the feedback law (transformed back to the original coordinates) (3.26), (3.27) with

$$v_1(t) = \frac{63\pi}{256}\left(30\left((\sigma + 4\pi^2)c_1(t) + f_1(t)\right) + 11\left((\sigma + \pi^2)c_2(t) + f_2(t)\right)\right)$$
$$v_2(t) = -\frac{495\pi}{256}\left(14\left((\sigma + 4\pi^2)c_1(t) + f_1(t)\right) + 3\left((\sigma + \pi^2)c_2(t) + f_2(t)\right)\right) \quad (3.30)$$

where $\sigma > 0$ is an (arbitrary) constant,

$$f_n(t) := \int_0^1 \sin(n\pi x)F(u(t,x))dx, \text{ for } n = 1, 2 \quad (3.31)$$

$$c_n(t) = \int_0^1 \sin(n\pi x)u(t,x)dx - \frac{4(-1)^n n}{(25-4n^2)\pi}y_1(t) - \frac{4(-1)^n n}{(49-4n^2)\pi}y_2(t), \text{ for } n=1,2 \qquad (3.32)$$

achieves global exponential stabilization of system (3.22), provided that the inequality

$$\bar{L} < 0.299 \qquad (3.33)$$

holds. ◁

## 4. Proofs of Main Results

We start with the proof of Lemma 2.1.

**Proof of Lemma 2.1:** We first show that $B_{n,i} = -\langle \varphi_i, \phi_n \rangle \neq 0$ for $n=1,...,N$, $i=1,...,j$. Indeed, using integration by parts, (2.1), (2.9), (2.10) and the fact that $A\phi_n = \lambda_n \phi_n$, we get for $i=1,...,j$:

$$\begin{aligned}\mu_i \langle \varphi_i, \phi_n \rangle &= \langle \mu_i \varphi_i, \phi_n \rangle = -\int_0^1 (p(x)\varphi_i'(x))' \phi_n(x)dx + \int_0^1 q(x)\varphi_i(x)\phi_n(x)dx \\ &= -p(1)(\varphi_i'(1)\phi_n(1) - \varphi_i(1)\phi_n'(1)) + p(0)(\varphi_i'(0)\phi_n(0) - \varphi_i(0)\phi_n'(0)) + \lambda_n \langle \varphi_i, \phi_n \rangle\end{aligned} \qquad (4.1)$$

Since the system $x_1\varphi_i(0) + x_2\varphi_i'(0) = 0$, $x_1\phi_n(0) + x_2\phi_n'(0) = 0$ admits a non-zero solution (namely, $(x_1, x_2) = (b_1, b_2)$; recall (2.2), (2.10)), it follows that $\varphi_i'(0)\phi_n(0) - \varphi_i(0)\phi_n'(0) = 0$. Moreover, since $a_1\phi_n(1) + a_2\phi_n'(1) = 0$, we obtain from (2.3) and (2.10) that $\varphi_i'(1)\phi_n(1) - \varphi_i(1)\phi_n'(1) = a_2\phi_n(1) - a_1\phi_n'(1)$. It follows from (4.1) that

$$B_{n,i} = -\langle \varphi_i, \phi_n \rangle = p(1)\frac{a_2\phi_n(1) - a_1\phi_n'(1)}{\mu_i - \lambda_n} \qquad (4.2)$$

The fact that $B_{n,i} \neq 0$ is shown by means of a contradiction. If $B_{n,i} = 0$ then we obtain $a_1\phi_n(1) + a_2\phi_n'(1) = 0$ and $a_2\phi_n(1) - a_1\phi_n'(1) = 0$. These two equations imply that $\phi_n(1) = \phi_n'(1) = 0$, which gives that $\phi_n(x) \equiv 0$ (uniqueness of solution for the initial-value problem $A\phi_n = \lambda_n\phi_n$ with $\phi_n(1) = \phi_n'(1) = 0$); a contradiction.

We are now ready to turn to the proof of the lemma. It suffices to show that the single-input, linear, time-invariant system

$$\begin{aligned}\dot{z} &= Cz + B_1 v_1 \\ z &\in \Re^N, v_1 \in \Re\end{aligned} \qquad (4.3)$$

is controllable. The Kalman rank controllability test for system (4.3) gives the square matrix

$$Q = \begin{bmatrix} B_{1,1} & (-\lambda_1)B_{1,1} & \cdots & (-\lambda_1)^{N-1}B_{1,1} \\ B_{2,1} & (-\lambda_2)B_{2,1} & \cdots & (-\lambda_2)^{N-1}B_{2,1} \\ \vdots & \vdots & & \vdots \\ B_{N,1} & (-\lambda_N)B_{N,1} & \cdots & (-\lambda_N)^{N-1}B_{N,1} \end{bmatrix} = diag(B_{1,1},...,B_{N,1})\begin{bmatrix} 1 & (-\lambda_1) & \cdots & (-\lambda_1)^{N-1} \\ 1 & (-\lambda_2) & \cdots & (-\lambda_2)^{N-1} \\ \vdots & \vdots & & \vdots \\ 1 & (-\lambda_N) & \cdots & (-\lambda_N)^{N-1} \end{bmatrix} \qquad (4.4)$$

Since $B_{n,1} = -\langle \varphi_1, \phi_n \rangle \neq 0$ for $n = 1, \ldots, N$, it follows that the matrix $diag(B_{1,1}, \ldots, B_{N,1})$ is invertible with non-zero determinant. Since the matrix
$$\begin{bmatrix} 1 & (-\lambda_1) & \cdots & (-\lambda_1)^{N-1} \\ 1 & (-\lambda_2) & \cdots & (-\lambda_2)^{N-1} \\ \vdots & \vdots & & \vdots \\ 1 & (-\lambda_N) & \cdots & (-\lambda_N)^{N-1} \end{bmatrix}$$
is a Vandermonde matrix and since $\lambda_1 < \lambda_2 < \ldots < \lambda_n < \ldots$, it follows that the matrix
$$\begin{bmatrix} 1 & (-\lambda_1) & \cdots & (-\lambda_1)^{N-1} \\ 1 & (-\lambda_2) & \cdots & (-\lambda_2)^{N-1} \\ \vdots & \vdots & & \vdots \\ 1 & (-\lambda_N) & \cdots & (-\lambda_N)^{N-1} \end{bmatrix}$$
has a non-zero determinant. It follows from (4.4) that the determinant of the matrix $Q$ is non-zero. Thus $rank(Q) = N$ and system (4.3) is controllable. The proof is complete. ◁

We next provide the proof of Theorem 2.2.

**Proof of Theorem 2.2:** Using Parseval's identity, the fact that $\lambda_1 < \lambda_2 < \ldots < \lambda_n < \ldots$, equation (2.16) and definition (2.25), we obtain for all $(w, v) \in D \times \Re^j$:

$$\|w - Pw\|^2 = \|w\|^2 - \sum_{n=1}^{N} \langle \phi_n, w \rangle^2 = \sum_{n=N+1}^{\infty} \langle \phi_n, w \rangle^2$$
$$\|Pw\|^2 = \sum_{n=1}^{N} \langle \phi_n, w \rangle^2$$
(4.5)

$$\left\langle w - Pw, -Aw + PAw - \sum_{i=1}^{j} v_i \varphi_i + \sum_{i=1}^{j} v_i P\varphi_i \right\rangle = -\sum_{n=N+1}^{\infty} \lambda_n \langle \phi_n, w \rangle^2 - \sum_{i=1}^{j} v_i \sum_{n=N+1}^{\infty} \langle \phi_n, w \rangle \langle \phi_n, \varphi_i \rangle$$
$$\leq -\lambda_{N+1}\left(\|w\|_2^2 - \sum_{n=1}^{N} \langle \phi_n, w \rangle^2\right) - \sum_{n=N+2}^{\infty} (\lambda_n - \lambda_{N+1}) \langle \phi_n, w \rangle^2 - \sum_{i=1}^{j} v_i \langle w - Pw, \varphi_i \rangle$$
(4.6)

Therefore, we obtain from (2.24), (2.29) and (4.6) for all $(w, y, v) \in D \times \Re^j \times \Re^j$:

$$\dot{V} \leq -\sigma \sum_{n=1}^{N} \langle \phi_n, w \rangle^2 - \sum_{i=1}^{j}\left(v_i - \sum_{n=1}^{N} K_{i,n} \langle \phi_n, w \rangle\right)\langle Gw, \varphi_i \rangle - \sum_{i=1}^{j} \omega_i \mu_i y_i^2 + \sum_{i=1}^{j} \omega_i y_i v_i$$
$$-\gamma \lambda_{N+1} \|w - Pw\|^2 - \gamma \sum_{i=1}^{j} v_i \langle w - Pw, \varphi_i \rangle - \gamma \sum_{n=N+2}^{\infty} (\lambda_n - \lambda_{N+1}) \langle \phi_n, w \rangle^2$$
(4.7)

Rearranging the terms in the right hand side of (4.7), we obtain for all $(w, y, v) \in D \times \Re^j \times \Re^j$:

$$\dot{V} \leq -\sigma \sum_{n=1}^{N} \langle \phi_n, w \rangle^2 - \sum_{i=1}^{j}\left(\langle Gw, \varphi_i \rangle + \gamma \langle w - Pw, \varphi_i \rangle - \omega_i y_i\right)\left(v_i - \sum_{n=1}^{N} K_{i,n} \langle \phi_n, w \rangle\right)$$
$$-\gamma \lambda_{N+1} \|w - Pw\|^2 - \gamma \sum_{n=N+2}^{\infty} (\lambda_n - \lambda_{N+1}) \langle \phi_n, w \rangle^2 - \sum_{i=1}^{j} \omega_i \mu_i y_i^2$$
$$-\gamma \sum_{i=1}^{j}\left(\sum_{n=1}^{N} K_{i,n} \langle \phi_n, w \rangle\right)\langle w - Pw, \varphi_i \rangle + \sum_{i=1}^{j} \omega_i y_i \left(\sum_{n=1}^{N} K_{i,n} \langle \phi_n, w \rangle\right)$$
(4.8)

We next use the (Young) inequalities

$$\omega_i |y_i| |K_{i,n} \langle \phi_n, w \rangle| \leq \frac{\sigma}{4j} \langle \phi_n, w \rangle^2 + \frac{j}{\sigma} \omega_i^2 y_i^2 K_{i,n}^2$$

$$\gamma |K_{i,n} \langle \phi_n, w \rangle| |\langle w - Pw, \varphi_i \rangle| \leq \frac{\sigma}{4j} \langle \phi_n, w \rangle^2 + \frac{j}{\sigma} \gamma^2 K_{i,n}^2 \langle w - Pw, \varphi_i \rangle^2$$

in conjunction with (4.8) and we get the estimate for all $(w, y, v) \in D \times \mathfrak{R}^j \times \mathfrak{R}^j$:

$$\begin{aligned}\dot{V} \leq &-\frac{\sigma}{2} \sum_{n=1}^{N} \langle \phi_n, w \rangle^2 - \sum_{i=1}^{j} \left( \langle Gw, \varphi_i \rangle + \gamma \langle w - Pw, \varphi_i \rangle - \omega_i y_i \right) \left( v_i - \sum_{n=1}^{N} K_{i,n} \langle \phi_n, w \rangle \right) \\ &-\gamma \lambda_{N+1} \|w - Pw\|^2 - \gamma \sum_{n=N+2}^{\infty} (\lambda_n - \lambda_{N+1}) \langle \phi_n, w \rangle^2 - \sum_{i=1}^{j} \left( \omega_i \mu_i - \frac{j}{\sigma} \omega_i^2 |K_i|^2 \right) y_i^2 \\ &+ \frac{j}{\sigma} \gamma^2 \sum_{i=1}^{j} \langle w - Pw, \varphi_i \rangle^2 |K_i|^2 \end{aligned} \quad (4.9)$$

The Cauchy-Schwarz inequality $|\langle w - Pw, \varphi_i \rangle| \leq \|w - Pw\| \|\varphi_i\|$ combined with (4.9) gives the inequality:

$$\begin{aligned}\dot{V} \leq &-\frac{\sigma}{2} \sum_{n=1}^{N} \langle \phi_n, w \rangle^2 - \sum_{i=1}^{j} \left( \langle Gw, \varphi_i \rangle + \gamma \langle w - Pw, \varphi_i \rangle - \omega_i y_i \right) \left( v_i - \sum_{n=1}^{N} K_{i,n} \langle \phi_n, w \rangle \right) \\ &-\left( \gamma \lambda_{N+1} - \frac{j}{\sigma} \gamma^2 \sum_{i=1}^{j} \|\varphi_i\|^2 |K_i|^2 \right) \|w - Pw\|^2 - \gamma \sum_{n=N+2}^{\infty} (\lambda_n - \lambda_{N+1}) \langle \phi_n, w \rangle^2 \\ &-\sum_{i=1}^{j} \left( \omega_i \mu_i - \frac{j}{\sigma} \omega_i^2 |K_i|^2 \right) y_i^2 \end{aligned} \quad (4.10)$$

It should be noticed that inequalities (2.30), (2.31) imply the inequalities $\gamma \lambda_{N+1} - \frac{j}{\sigma} \gamma^2 \sum_{i=1}^{j} \|\varphi_i\|^2 |K_i|^2 \geq \frac{1}{2} \gamma \lambda_{N+1}$ and $\omega_i \mu_i - \frac{j}{\sigma} \omega_i^2 |K_i|^2 \geq \frac{1}{2} \omega_i \mu_i$ for $i = 1, ..., j$. Consequently, we get from (4.10) for all $(w, y, v) \in D \times \mathfrak{R}^j \times \mathfrak{R}^j$:

$$\begin{aligned}\dot{V} \leq &-\frac{\sigma}{2} \sum_{n=1}^{N} \langle \phi_n, w \rangle^2 - \sum_{i=1}^{j} \left( \langle Gw, \varphi_i \rangle + \gamma \langle w - Pw, \varphi_i \rangle - \omega_i y_i \right) \left( v_i - \sum_{n=1}^{N} K_{i,n} \langle \phi_n, w \rangle \right) \\ &-\frac{\gamma \lambda_{N+1}}{2} \|w - Pw\|^2 - \gamma \sum_{n=N+2}^{\infty} (\lambda_n - \lambda_{N+1}) \langle \phi_n, w \rangle^2 - \sum_{i=1}^{j} \frac{\omega_i \mu_i}{2} y_i^2 \end{aligned} \quad (4.11)$$

The feedback law (2.33), (2.34) implies that the following equations hold for $i = 1, ..., j$:

$$v_i = \sum_{n=1}^{N} K_{i,n} \langle \phi_n, w \rangle + L_i \left( \langle Gw, \varphi_i \rangle + \gamma \langle w - Pw, \tilde{\varphi}_i \rangle - \omega_i y_i \right) \quad (4.12)$$

where

$$\tilde{\varphi}_i = \sum_{n=1}^{M} \langle \phi_n, \varphi_i \rangle \phi_n, \text{ for } i = 1, ..., j \quad (4.13)$$

Inserting (4.12) into (4.11) we obtain the following estimate for all $(w, y) \in D \times \mathfrak{R}^j$:

$$\dot{V} \leq -\frac{\sigma}{2}\sum_{n=1}^{N}\langle\phi_n,w\rangle^2 - \sum_{i=1}^{j}L_i\left(\langle Gw,\varphi_i\rangle + \gamma\langle w-Pw,\tilde{\varphi}_i\rangle - \omega_i y_i\right)^2$$
$$-\frac{\gamma\lambda_{N+1}}{2}\|w-Pw\|^2 - \gamma\sum_{n=N+2}^{\infty}(\lambda_n-\lambda_{N+1})\langle\phi_n,w\rangle^2 - \sum_{i=1}^{j}\frac{\omega_i\mu_i}{2}y_i^2 \qquad (4.14)$$
$$+\sum_{i=1}^{j}L_i\gamma\left|\langle w-Pw,\varphi_i-\tilde{\varphi}_i\rangle\right|\left|\langle Gw,\varphi_i\rangle + \gamma\langle w-Pw,\tilde{\varphi}_i\rangle - \omega_i y_i\right|$$

Using the (Young) inequality

$$\gamma\left|\langle w-Pw,\varphi_i-\tilde{\varphi}_i\rangle\right|\left|\langle Gw,\varphi_i\rangle + \gamma\langle w-Pw,\tilde{\varphi}_i\rangle - \omega_i y_i\right|$$
$$\leq\left|\langle Gw,\varphi_i\rangle + \gamma\langle w-Pw,\tilde{\varphi}_i\rangle - \omega_i y_i\right|^2 + \frac{\gamma^2}{4}\left|\langle w-Pw,\varphi_i-\tilde{\varphi}_i\rangle\right|^2$$

and the fact $M\geq N+1$, we obtain from (4.14) for all $(w,y)\in D\times\Re^j$:

$$\dot{V}\leq -\frac{\sigma}{2}\sum_{n=1}^{N}\langle\phi_n,w\rangle^2 - \frac{\gamma\lambda_{N+1}}{2}\|w-Pw\|^2 - \gamma\sum_{n=M+1}^{\infty}(\lambda_n-\lambda_{N+1})\langle\phi_n,w\rangle^2$$
$$-\sum_{i=1}^{j}\frac{\omega_i\mu_i}{2}y_i^2 + \sum_{i=1}^{j}L_i\frac{\gamma^2}{4}\left|\langle w-Pw,\varphi_i-\tilde{\varphi}_i\rangle\right|^2 \qquad (4.15)$$

Using the fact that $\langle w-Pw,\varphi_i-\tilde{\varphi}_i\rangle = \left\langle w-\sum_{n=1}^{M}\phi_n\langle w,\phi_n\rangle,\varphi_i-\tilde{\varphi}_i\right\rangle$ (a consequence of definition (4.13) and the fact that $M\geq N+1$) in conjunction with the Cauchy-Schwarz inequality and Parseval's identity, which gives

$$\left|\langle w-Pw,\varphi_i-\tilde{\varphi}_i\rangle\right|\leq\left(\sum_{n=M+1}^{\infty}\langle\phi_n,w\rangle^2\right)^{1/2}\left(\sum_{n=M+1}^{\infty}\langle\phi_n,\varphi_i\rangle^2\right)^{1/2}$$

we obtain from (4.15) the following estimate for all $(w,y)\in D\times\Re^j$:

$$\dot{V}\leq -\frac{\sigma}{2}\sum_{n=1}^{N}\langle\phi_n,w\rangle^2 - \frac{\gamma\lambda_{N+1}}{2}\|w-Pw\|^2 - \sum_{i=1}^{j}\frac{\omega_i\mu_i}{2}y_i^2$$
$$-\gamma\sum_{n=M+1}^{\infty}\left(\lambda_n-\lambda_{N+1}-\frac{\gamma}{4}\sum_{i=1}^{j}L_i\sum_{l=M+1}^{\infty}\langle\phi_l,\varphi_i\rangle^2\right)\langle\phi_n,w\rangle^2 \qquad (4.16)$$

The facts that $\lambda_1<\lambda_2<\ldots<\lambda_n<\ldots$, $4(\lambda_{M+1}-\lambda_{N+1})\geq\gamma\sum_{i=1}^{j}L_i\sum_{n=M+1}^{\infty}\langle\phi_n,\varphi_i\rangle^2$ in conjunction with (4.16) and (4.5) imply that inequality (2.32) holds.

Finally, let arbitrary $w_0\in D$, $y_0\in\Re^j$ be given. Existence-uniqueness of the solution $w\in C^0(\Re_+\times[0,1])\cap C^1((0,+\infty)\times[0,1])$, $y\in C^1(\Re_+;\Re^j)$ with $w[t]\in D\cap C^2([0,1])$ for all $t>0$, of the initial-boundary value problem (2.13), (2.14), (2.15), (2.33) with initial condition $w[0]=w_0$, $y(0)=y_0$ follows from Theorem 4.5 on page 73 of the book [9] and the fact that the functions $k_i$, $i=1,\ldots,j$ are of class $D$ (recall (2.34)). Therefore, by virtue of (2.32) the function $h(t)=V(w[t],y(t))$ is of class $C^0(\Re_+)\cap C^1((0,+\infty))$ and satisfies the following differential inequality for all $t>0$:

$$\dot{h}(t)\leq -\frac{1}{2}\sum_{i=1}^{j}\omega_i\mu_i y_i^2(t) - \frac{\min(\gamma\lambda_{N+1},\sigma)}{2}\|w[t]\|^2 \qquad (4.17)$$

Using (2.27), we conclude from (4.17) that there exists a constant $\bar{\sigma} > 0$ (independent of the particular solution $(w, y)$) such that the following differential inequality holds for all $t > 0$:

$$\dot{h}(t) \leq -2\bar{\sigma} h(t) \tag{4.18}$$

The differential inequality (4.18) in conjunction with continuity of $h(t) = V(w[t], y(t))$ at $t = 0$ implies the following estimate:

$$V(w[t], y(t)) \leq \exp(-2\bar{\sigma} t) V(w_0, y_0), \text{ for all } t \geq 0 \tag{4.19}$$

The exponential stability estimate (2.35) for certain appropriate constant $\bar{K} > 0$ (independent of the particular solution $(w, y)$) is a direct consequence of (4.19). The proof is complete. ◁

**Proof of Theorem 3.1:** Let arbitrary $(w, y) \in D \times \Re^N$ be given. In what follows we will use the notation $u = w + \sum_{i=1}^{N} \varphi_i y_i$. Define

$$h(\zeta) := \frac{2\zeta}{(1-\zeta)(1+\bar{L}^2)(1+\kappa N)}, \text{ for } \zeta \in (-\infty, 1) \tag{4.20}$$

and notice that by virtue of (3.13), (3.14) there exists sufficiently small $\zeta \in (0, 1)$ such that

$$\frac{\mu_i^2}{h(\zeta) \sum_{m=1}^{N} g_{i,m}^2 (\sigma - \lambda_m)^2 + 2 \sum_{m=1}^{N} g_{i,m}^2} - N\bar{L}^2 \left(1 + \frac{1}{\kappa}\right) \|\varphi_i\|^2 > 0, \text{ for } i = 1, \ldots, N \tag{4.21}$$

$$\frac{\lambda_{N+1}^2}{1 + Nh(\zeta) \sum_{i=1}^{N} \sum_{m=1}^{N} \|\varphi_i\|^2 g_{i,m}^2 (\sigma - \lambda_m)^2 + 2N \sum_{i=1}^{N} \sum_{m=1}^{N} \|\varphi_i\|^2 g_{i,m}^2} - \bar{L}^2 (1 + \kappa N) > 0 \tag{4.22}$$

Define:

$$\beta := N^{-1}(1-\zeta)\sigma R \tag{4.23}$$

$$\varepsilon := \frac{1}{2N\zeta}, \quad \gamma := \frac{\lambda_{N+1}}{\frac{1}{2N\zeta} + \sum_{i=1}^{N} \sum_{m=1}^{N} \|\varphi_i\|^2 \left(\frac{g_{i,m}^2 (\sigma - \lambda_m)^2}{\beta} + \frac{g_{i,m}^2}{\zeta}\right)}, \quad R := \sigma^{-1} N(1 + \bar{L}^2)(1 + \kappa N) \tag{4.24}$$

$$\omega_i = \frac{\mu_i}{\sum_{m=1}^{N} \left(\frac{g_{i,m}^2 (\sigma - \lambda_m)^2}{\beta} + \frac{g_{i,m}^2}{\zeta}\right)}, \text{ for } i = 1, \ldots, N \tag{4.25}$$

Using the Cauchy-Schwarz inequality, (3.9), (3.10), (4.6) and (3.17), we get:

$$\begin{aligned}
\dot{V} \leq & -\sigma R \sum_{m=1}^{N} \langle \phi_m, w \rangle^2 - \gamma \lambda_{N+1} \|w - Pw\|^2 + \gamma \|w - Pw\| \|F(u) - PF(u)\| \\
& -\gamma \sum_{i=1}^{N} \sum_{m=1}^{N} g_{i,m} (\sigma - \lambda_m) \langle \phi_m, w \rangle \langle w - Pw, \varphi_i \rangle - \gamma \sum_{i=1}^{N} \sum_{m=1}^{N} g_{i,m} \langle \phi_m, F(u) \rangle \langle w - Pw, \varphi_i \rangle \\
& - \sum_{i=1}^{N} \omega_i \mu_i y_i^2 + \sum_{i=1}^{N} \sum_{m=1}^{N} g_{i,m} (\sigma - \lambda_m) \omega_i y_i \langle \phi_m, w \rangle + \sum_{i=1}^{N} \sum_{m=1}^{N} g_{i,m} \omega_i y_i \langle \phi_m, F(u) \rangle
\end{aligned} \tag{4.26}$$

Using the inequalities

$$\gamma \|w - Pw\| \|F(u) - PF(u)\| \le \frac{\varepsilon}{2}\gamma^2 \|w - Pw\|^2 + \frac{1}{2\varepsilon}\|F(u) - PF(u)\|^2$$

$$\left|\gamma g_{i,m}(\sigma - \lambda_m)\langle \phi_m, w\rangle \langle w - Pw, \varphi_i\rangle\right| \le \frac{\beta}{2}\langle \phi_m, w\rangle^2 + \frac{\gamma^2}{2\beta} g_{i,m}^2 (\sigma - \lambda_m)^2 \langle w - Pw, \varphi_i\rangle^2$$

$$\left|\gamma g_{i,m}\langle \phi_m, F(u)\rangle \langle w - Pw, \varphi_i\rangle\right| \le \frac{\zeta}{2}\langle \phi_m, F(u)\rangle^2 + \frac{\gamma^2}{2\zeta} g_{i,m}^2 \langle w - Pw, \varphi_i\rangle^2$$

$$\left|g_{i,m}(\sigma - \lambda_m)\langle \phi_m, w\rangle \omega_i y_i\right| \le \frac{\beta}{2}\langle \phi_m, w\rangle^2 + \frac{\omega_i^2}{2\beta} g_{i,m}^2 (\sigma - \lambda_m)^2 y_i^2$$

$$\left|g_{i,m}\langle \phi_m, F(u)\rangle \omega_i y_i\right| \le \frac{\zeta}{2}\langle \phi_m, F(u)\rangle^2 + \frac{\omega_i^2}{2\zeta} g_{i,m}^2 y_i^2$$

the Cauchy-Schwarz inequality $\left|\langle w - Pw, \varphi_i\rangle\right| \le \|w - Pw\| \|\varphi_i\|$ and the fact that $\|PF(u)\|^2 = \sum_{m=1}^{N}\langle \phi_m, F(u)\rangle^2$, we get from (4.26) the following inequality:

$$\begin{aligned}
\dot{V} \le &-\sum_{m=1}^{N}(\sigma R - N\beta)\langle \phi_m, w\rangle^2 - \sum_{i=1}^{N}\left(\omega_i \mu_i - \frac{\omega_i^2}{2}\sum_{m=1}^{N}\left(\frac{g_{i,m}^2(\sigma - \lambda_m)^2}{\beta} + \frac{g_{i,m}^2}{\zeta}\right)\right)y_i^2 \\
&-\left(\gamma \lambda_{N+1} - \frac{1}{2}\gamma^2\left(\varepsilon + \sum_{i=1}^{N}\sum_{m=1}^{N}\|\varphi_i\|_2^2\left(\frac{g_{i,m}^2(\sigma - \lambda_m)^2}{\beta} + \frac{g_{i,m}^2}{\zeta}\right)\right)\right)\|w - Pw\|^2 \\
&+ N\zeta \|PF(u)\|^2 + \frac{1}{2\varepsilon}\|F(u) - PF(u)\|^2
\end{aligned} \tag{4.27}$$

Since $N\zeta = \frac{1}{2\varepsilon}$ (a consequence of (4.24)) and since $\|PF(u)\|^2 + \|F(u) - PF(u)\|^2 = \|F(u)\|^2$ (a consequence of definition (2.25)), we obtain from (3.4) and (4.27):

$$\begin{aligned}
\dot{V} \le &-\sum_{m=1}^{N}(\sigma R - N\beta)\langle \phi_m, w\rangle^2 - \sum_{i=1}^{N}\left(\omega_i \mu_i - \frac{\omega_i^2}{2}\sum_{m=1}^{N}\left(\frac{g_{i,m}^2(\sigma - \lambda_m)^2}{\beta} + \frac{g_{i,m}^2}{\zeta}\right)\right)y_i^2 \\
&-\left(\gamma \lambda_{N+1} - \frac{1}{2}\gamma^2\left(\frac{1}{2N\zeta} + \sum_{i=1}^{N}\sum_{m=1}^{N}\|\varphi_i\|^2\left(\frac{g_{i,m}^2(\sigma - \lambda_m)^2}{\beta} + \frac{g_{i,m}^2}{\zeta}\right)\right)\right)\|w - Pw\|^2 + N\zeta \bar{L}^2 \|u\|^2
\end{aligned} \tag{4.28}$$

Since $u = w + \sum_{i=1}^{N}\varphi_i y_i$ and since $\langle \varphi_i, \varphi_m\rangle = 0$ for $i \ne m$ (recall (3.8)), it follows that $\|u\|^2 = \|w\|^2 + 2\sum_{i=1}^{N}\langle \varphi_i, w\rangle y_i + \sum_{i=1}^{N}\|\varphi_i\|^2 y_i^2$. Using the Cauchy-Schwarz inequality $\left|\langle \varphi_i, w\rangle\right| \le \|\varphi_i\| \|w\|$ and the fact that $\|\varphi_i\| \|w\| |y_i| \le \frac{1}{2\kappa}\|\varphi_i\|^2 y_i^2 + \frac{\kappa}{2}\|w\|^2$, we obtain the estimate $\|u\|^2 \le (1 + \kappa N)\|w\|^2 + \left(1 + \frac{1}{\kappa}\right)\sum_{i=1}^{N}\|\varphi_i\|^2 y_i^2$. Therefore, we obtain from (4.5) and (4.28):

$$\dot{V} \leq -\sum_{m=1}^{N}\left(\sigma R - N\beta - N\zeta \bar{L}^2(1+\kappa N)\right)\langle \phi_m, w\rangle^2$$

$$-\sum_{i=1}^{N}\left(\omega_i \mu_i - \frac{\omega_i^2}{2}\sum_{m=1}^{N}\left(\frac{g_{i,m}^2(\sigma-\lambda_m)^2}{\beta} + \frac{g_{i,m}^2}{\zeta}\right) - N\zeta\bar{L}^2\left(1+\frac{1}{\kappa}\right)\|\varphi_i\|^2\right)y_i^2 \quad (4.29)$$

$$-\left(\gamma\lambda_{N+1} - \frac{1}{2}\gamma^2\left(\frac{1}{2N\zeta} + \sum_{i=1}^{N}\sum_{m=1}^{N}\|\varphi_i\|^2\left(\frac{g_{i,m}^2(\sigma-\lambda_m)^2}{\beta} + \frac{g_{i,m}^2}{\zeta}\right)\right) - N\zeta\bar{L}^2(1+\kappa N)\right)\|w - Pw\|^2$$

Substituting definition (4.24) of $\gamma$ and definitions (4.25) of $\omega_i$ into (4.29), we obtain the inequality:

$$\dot{V} \leq -\left(\sigma R - N\beta - N\zeta\bar{L}^2(1+\kappa N)\right)\sum_{m=1}^{N}\langle\phi_m, w\rangle^2$$

$$-\zeta\sum_{i=1}^{N}\left(\frac{\mu_i^2}{2\sum_{m=1}^{N}\left(\frac{\zeta g_{i,m}^2(\sigma-\lambda_m)^2}{\beta} + g_{i,m}^2\right)} - N\bar{L}^2\left(1+\frac{1}{\kappa}\right)\|\varphi_i\|^2\right)y_i^2 \quad (4.30)$$

$$-N\zeta\left(\frac{\lambda_{N+1}^2}{1+2N\sum_{i=1}^{N}\sum_{m=1}^{N}\|\varphi_i\|^2\left(\frac{\zeta g_{i,m}^2(\sigma-\lambda_m)^2}{\beta} + g_{i,m}^2\right)} - \bar{L}^2(1+\kappa N)\right)\|w-Pw\|^2$$

Using (4.20) and substituting definition (4.23) of $\beta$ and definitions (4.24) of $R$ into (4.30), we obtain the inequality:

$$\dot{V} \leq -\zeta N(1+\kappa N)\sum_{m=1}^{N}\langle\phi_m, w\rangle^2 - \zeta\sum_{i=1}^{N}\left(\frac{\mu_i^2}{h(\zeta)\sum_{m=1}^{N}g_{i,m}^2(\sigma-\lambda_m)^2 + 2\sum_{m=1}^{N}g_{i,m}^2} - N\bar{L}^2\left(1+\frac{1}{\kappa}\right)\|\varphi_i\|^2\right)y_i^2$$

$$-N\zeta\left(\frac{\lambda_{N+1}^2}{1+Nh(\zeta)\sum_{i=1}^{N}\sum_{m=1}^{N}\|\varphi_i\|^2 g_{i,m}^2(\sigma-\lambda_m)^2 + 2N\sum_{i=1}^{N}\sum_{m=1}^{N}\|\varphi_i\|^2 g_{i,m}^2} - \bar{L}^2(1+\kappa N)\right)\|w-Pw\|^2 \quad (4.31)$$

Inequality (3.16) for certain appropriate constant $\theta > 0$ is a direct consequence of inequalities (4.21), (4.22), (4.31) and equation (4.5).

From this point we follow the same procedure as in the proof of Theorem 2.2 and we derive the exponential stability estimate (2.35). The proof is complete. ◁

**Proof of Theorem 3.2:** Let arbitrary $(w, y) \in D \times \Re^N$ be given. In what follows we will use the notation $u = w + \sum_{i=1}^{N}\varphi_i y_i$. Notice that by virtue of (3.18), (3.19), (3.20) there exists sufficiently small $a \in (0,1)$ such that

$$\sigma^2 - \bar{L}^2(1+\kappa N) > a \quad (4.32)$$

$$\mu_i^2 > \left(1+\frac{1}{\kappa}\right)\frac{2N\overline{L}^2\|\varphi_i\|^2}{\sigma^2 - a - \overline{L}^2(1+\kappa N)}\left(\varepsilon + \sum_{m=1}^{N} g_{i,m}^2(\sigma - \lambda_m)^2\right), \text{ for } i = 1,...,N \tag{4.33}$$

$$\lambda_{N+1}^2 > \overline{L}^2(1+\kappa N)\left(1 + 2N\sum_{i=1}^{N}\sum_{m=1}^{N}\frac{\|\varphi_i\|^2 g_{i,m}^2(\sigma - \lambda_m)^2}{\sigma^2 - a - \overline{L}^2(1+\kappa N)}\right) \tag{4.34}$$

Define:

$$\beta := \frac{\sigma^2 - a - \overline{L}^2(1+\kappa N)}{2N} \tag{4.35}$$

$$\gamma := \frac{\beta \lambda_{N+1}}{\beta + \sum_{i=1}^{N}\sum_{m=1}^{N}\|\varphi_i\|^2 g_{i,m}^2(\sigma - \lambda_m)^2}, \quad R := \sigma \tag{4.36}$$

$$\omega_i := \frac{\beta \mu_i}{\varepsilon + \sum_{m=1}^{N} g_{i,m}^2(\sigma - \lambda_m)^2}, \text{ for } i = 1,...,N \tag{4.37}$$

Notice that (4.32) implies that $\beta$ as defined by (4.35) is positive. Using (3.11), (3.12), (4.6) and (3.17), we obtain:

$$\begin{aligned}\dot{V} \leq &-\sum_{m=1}^{N}\sigma R\langle\phi_m,w\rangle^2 + \sum_{m=1}^{N}R\langle\phi_m,w\rangle\langle\phi_m,F(u)\rangle - \gamma\lambda_{N+1}\|w - Pw\|^2 \\ &+ \gamma\langle w - Pw, F(u) - PF(u)\rangle + \sum_{i=1}^{N}\sum_{m=1}^{N}g_{i,m}(\sigma - \lambda_m)\omega_i y_i\langle\phi_m,w\rangle \\ &- \gamma\sum_{i=1}^{N}\sum_{m=1}^{N}g_{i,m}(\sigma - \lambda_m)\langle\phi_m,w\rangle\langle w - Pw, \varphi_i\rangle - \sum_{i=1}^{N}\omega_i\mu_i y_i^2\end{aligned} \tag{4.38}$$

Using the inequalities

$$\gamma\|w - Pw\|\|F(u) - PF(u)\| \leq \frac{1}{2}\gamma^2\|w - Pw\|^2 + \frac{1}{2}\|F(u) - PF(u)\|^2$$

$$\left|R\langle\phi_m,w\rangle\langle\phi_m,F(u)\rangle\right| \leq \frac{1}{2}R^2\langle\phi_m,w\rangle^2 + \frac{1}{2}\langle\phi_m,F(u)\rangle^2$$

$$\left|\gamma g_{i,m}(\sigma - \lambda_m)\langle\phi_m,w\rangle\langle w - Pw,\varphi_i\rangle\right| \leq \frac{\beta}{2}\langle\phi_m,w\rangle^2 + \frac{\gamma^2}{2\beta}g_{i,m}^2(\sigma - \lambda_m)^2\langle w - Pw,\varphi_i\rangle^2$$

$$\left|g_{i,m}(\sigma - \lambda_m)\langle\phi_m,w\rangle\omega_i y_i\right| \leq \frac{\beta}{2}\langle\phi_m,w\rangle^2 + \frac{\omega_i^2}{2\beta}g_{i,m}^2(\sigma - \lambda_m)^2 y_i^2$$

the Cauchy-Schwarz inequalities $\left|\langle w - Pw, F(u) - PF(u)\rangle\right| \leq \|w - Pw\|\|F(u) - PF(u)\|$, $\left|\langle w - Pw, \varphi_i\rangle\right| \leq \|w - Pw\|\|\varphi_i\|$ and the fact that $\|PF(u)\|^2 = \sum_{m=1}^{N}\langle\phi_m,F(u)\rangle^2$, we get from (4.38) the following inequality:

$$\dot{V} \leq -\sum_{m=1}^{N}\left(\sigma R-\frac{1}{2}R^{2}-N\beta\right)\langle\phi_{m},w\rangle^{2}+\frac{1}{2}\|PF(u)\|^{2}$$
$$-\left(\gamma\lambda_{N+1}-\frac{1}{2}\gamma^{2}-\frac{1}{2\beta}\gamma^{2}\sum_{i=1}^{N}\sum_{m=1}^{N}\|\varphi_{i}\|^{2}\,g_{i,m}^{2}(\sigma-\lambda_{m})^{2}\right)\|w-Pw\|^{2} \qquad (4.39)$$
$$+\frac{1}{2}\|F(u)-PF(u)\|^{2}-\sum_{i=1}^{N}\left(\omega_{i}\mu_{i}-\frac{\omega_{i}^{2}}{2\beta}\sum_{m=1}^{N}g_{i,m}^{2}(\sigma-\lambda_{m})^{2}\right)y_{i}^{2}$$

Since $\|PF(u)\|^{2}+\|F(u)-PF(u)\|^{2}=\|F(u)\|^{2}$ (a consequence of definition (2.25)), we obtain from (3.4) and (4.39):

$$\dot{V} \leq -\sum_{m=1}^{N}\left(\sigma R-\frac{1}{2}R^{2}-N\beta\right)\langle\phi_{m},w\rangle^{2}+\frac{\overline{L}^{2}}{2}\|u\|^{2}$$
$$-\left(\gamma\lambda_{N+1}-\frac{1}{2}\gamma^{2}-\frac{1}{2\beta}\gamma^{2}\sum_{i=1}^{N}\sum_{m=1}^{N}\|\varphi_{i}\|^{2}\,g_{i,m}^{2}(\sigma-\lambda_{m})^{2}\right)\|w-Pw\|^{2} \qquad (4.40)$$
$$-\sum_{i=1}^{N}\left(\omega_{i}\mu_{i}-\frac{\omega_{i}^{2}}{2\beta}\sum_{m=1}^{N}g_{i,m}^{2}(\sigma-\lambda_{m})^{2}\right)y_{i}^{2}$$

Since $u=w+\sum_{i=1}^{N}\varphi_{i}y_{i}$ and since $\langle\varphi_{i},\varphi_{m}\rangle=0$ for $i\neq m$ (recall (3.8)), it follows that $\|u\|^{2}=\|w\|^{2}+2\sum_{i=1}^{N}\langle\varphi_{i},w\rangle y_{i}+\sum_{i=1}^{N}\|\varphi_{i}\|^{2}\,y_{i}^{2}$. Using the Cauchy-Schwarz inequality $|\langle\varphi_{i},w\rangle|\leq\|\varphi_{i}\|\|w\|$ and the fact that $\|\varphi_{i}\|\|w\||y_{i}|\leq\frac{1}{2\kappa}\|\varphi_{i}\|^{2}\,y_{i}^{2}+\frac{\kappa}{2}\|w\|^{2}$, we obtain the estimate $\|u\|^{2}\leq(1+\kappa N)\|w\|^{2}+\left(1+\frac{1}{\kappa}\right)\sum_{i=1}^{N}\|\varphi_{i}\|^{2}\,y_{i}^{2}$. Therefore, we obtain from (4.5) and (4.40):

$$\dot{V} \leq -\sum_{m=1}^{N}\left(\sigma R-\frac{1}{2}R^{2}-N\beta-\frac{\overline{L}^{2}}{2}(1+\kappa N)\right)\langle\phi_{m},w\rangle^{2}$$
$$-\left(\gamma\lambda_{N+1}-\frac{1}{2}\gamma^{2}-\frac{1}{2\beta}\gamma^{2}\sum_{i=1}^{N}\sum_{m=1}^{N}\|\varphi_{i}\|^{2}\,g_{i,m}^{2}(\sigma-\lambda_{m})^{2}-\frac{\overline{L}^{2}}{2}(1+\kappa N)\right)\|w-Pw\|^{2} \qquad (4.41)$$
$$-\sum_{i=1}^{N}\left(\omega_{i}\mu_{i}-\frac{\omega_{i}^{2}}{2\beta}\left(\varepsilon+\sum_{m=1}^{N}g_{i,m}^{2}(\sigma-\lambda_{m})^{2}\right)-\left(1+\frac{1}{\kappa}\right)\frac{\overline{L}^{2}}{2}\|\varphi_{i}\|^{2}\right)y_{i}^{2}$$

Substituting the expressions for $\beta, R, \gamma, \omega_{1},...,\omega_{N}$ from (4.35), (4.36), (4.37) in (4.41), we obtain:

$$\dot{V} \leq -\frac{a}{2}\sum_{m=1}^{N}\langle\phi_{m},w\rangle^{2}-\frac{1}{2}\sum_{i=1}^{N}\left(\frac{\left(\sigma^{2}-a-\overline{L}^{2}(1+\kappa N)\right)\mu_{i}^{2}}{2\varepsilon N+2N\sum_{m=1}^{N}g_{i,m}^{2}(\sigma-\lambda_{m})^{2}}-\left(1+\frac{1}{\kappa}\right)\overline{L}^{2}\|\varphi_{i}\|^{2}\right)y_{i}^{2}$$
$$-\frac{1}{2}\left(\frac{\lambda_{N+1}^{2}}{1+2N\sum_{i=1}^{N}\sum_{m=1}^{N}\frac{\|\varphi_{i}\|^{2}\,g_{i,m}^{2}(\sigma-\lambda_{m})^{2}}{\sigma^{2}-a-\overline{L}^{2}(1+\kappa N)}}-\overline{L}^{2}(1+\kappa N)\right)\|w-Pw\|^{2} \qquad (4.42)$$

Inequality (3.16) for certain appropriate constant $\theta > 0$ is a direct consequence of inequalities (4.33), (4.34), (4.42) and equation (4.5).

From this point we follow the same procedure as in the proof of Theorem 2.2 and we derive the exponential stability estimate (2.35). The proof is complete. ◁

# 5. Concluding Remarks

The present paper showed that simple ("almost diagonal" with no double and triple integrals) CLFs for boundary controlled parabolic PDEs do exist. Moreover, these simple CLFs can give a novel methodology for boundary feedback design in parabolic PDEs. The proposed feedback design methodology is Lyapunov-based and the feedback controller is an "integral" controller with internal dynamics. It was also shown that the obtained simple CLFs can provide nonlinear boundary feedback laws which achieve global exponential stabilization of semilinear parabolic PDEs with nonlinearities that satisfy a linear growth condition.

Future research may address the construction of simple CLFs for boundary controlled systems of parabolic PDEs. Constructions of Lyapunov functionals for systems of parabolic PDEs were proposed in [15,16] and can be the basis of future work. Another research direction that could be explored is the further simplification of the CLF. It is possible that non-coercive CLFs (proposed in [17]) may provide ways to simplify a CLF for a boundary controlled parabolic PDE even further.

**Acknowledgments:** The author would like to thank Prof. M. Krstic for various discussions about the existence of "diagonal" CLFs for parabolic PDEs.